\documentclass{article}
\usepackage{amssymb}
\usepackage{amsmath}
\usepackage{setspace}

\newenvironment{Macaulay2}{
\begin{spacing}{0.7}
\begin{quote}
} {
\end{quote}
\end{spacing}
}

\newtheorem{theorem}{Theorem}

\newtheorem{lemma}[theorem]{Lemma}

\newtheorem{proposition}[theorem]{Proposition}
\newtheorem{remark}[theorem]{Remark}

\newenvironment{proof}[1][Proof]{\textbf{#1.} }{\ \rule{0.5em}{0.5em}\smallskip}

\newcommand{\p}{\partial}
\newcommand{\e}{\varepsilon}
\newcommand{\de}{\delta}
\newcommand{\len}{\mbox{length}} %length

\title{Algorithmic proofs of two theorems of Stafford}
\author{Anton Leykin}
\begin{document}
\maketitle
\begin{abstract}%%%%%%%%%%%%%%%%%%%%%%%%%%%%%%%%%%%%%%%%%%%%%%%%%%%%%%%%%%%%
Two classical results of Stafford say that every (left) ideal of the $n$-th
Weyl algebra $A_n$ can be generated by two elements,
and every holonomic $A_n$-module is cyclic, i.e. generated by one element.
We modify Stafford's original proofs to make the algorithmic computation
of these generators possible.
\end{abstract}%%%%%%%%%%%%%%%%%%%%%%%%%%%%%%%%%%%%%%%%%%%%%%%%%%%%%%%%%%%%%%

\section{Introduction}%%%%%%%%%%%%%%%%%%%%%%%%%%%%%%%%%%%%%%%%%%%%%%%%%%%%%%
Let $k$ is a field of characteristic $0$, and \(
A_n=A_{n}(k)=k\left\langle x_{1},...,x_{n},\partial
_{1},...,\partial _{n}\right\rangle  \) be the $n$-th Weyl
algebra, which is an associative \( k \)-algebra generated by \( x
\)'s and \( \partial \)'s with the relations \(\partial
_{i}x_{i}=x_{i}\partial _{i}+1 \) for all \( i \). This algebra
may be thought of as the algebra of linear differential operators
with polynomial coefficients.

There are several things that are nice about the Weyl algebra.
First of all the dimension theory can be developed for it; this is
done, for example, in Chapter 1 of Bj\"ork \cite{bjork}. It is
shown that the Gelfand-Kirillov dimension of $A_n$ equals $2n$,
moreover, if $M$ is a nontrivial $A_n$-module, then $n \leq
\dim{M} \leq 2n$. The modules of dimension $n$ (minimal possible
dimension) constitute the {\em Bernstein class}.

One of the distinctive properties of the modules in Bernstein
class, which are also called \emph{holonomic} modules, is their
finite length. Below we shall show that this property implies that
{\bf every holonomic module can be generated by one element}.

\smallskip
Another striking fact, which is very simple to state, but quite hard
to prove, is that {\bf for every left ideal of $A_n$ there exist
2 elements that generate it}.

\smallskip
Both statements were proved by Stafford in
\cite{stafford}, also these results appear in \cite{bjork}.
Unfortunately, the arguments given by Stafford can't be converted
to algorithms straightforwardly. There are several obstacles to
this, many of which one can overcome with the theory of Gr\"obner
bases for Weyl algebras. However, the main difficulty is that both
proofs contain an operation of taking an irreducible
submodule of an $A_n$-module. To our best knowledge, there doesn't
exist an algorithm for this; moreover, even if such algorithm is invented
one should expect it to be quite involved.

We were able to modify the original proofs in such a way that
computations are possible and implemented the corresponding
algorithms in the computer algebra system \emph{Macaulay 2} \cite{MTwo}.

We have to mention that in their recent paper \cite{schmale}
Hillebrand and Schmale construct another effective modification of
Stafford's proof which leads to an algorithm. We shall discuss the
differences of their and our approaches in the last section.

\section{Notation Table} %%%%%%%%%%%%%%%%%%%%%%%%%%%%%%%%%%%%%%%%%%%%%%%%
For the convenience of the reader we provide the notation lookup
table. All of the symbols listed below show up sooner or later in
the paper along with more detailed definitions.
\begin{eqnarray*}
k &\hbox{is}& \hbox{a (commutative) field of characteristic $0$},\\
A_r&=&A_r(k)=k\langle x_1,...,x_r,\p_1,...,\p_r\rangle, \\
A &\hbox{is}& \hbox{a simple ring of infinite length as
    a left module over itself}, \\
D &\hbox{is}& \hbox{a skew field of characteristic $0$},\\
K &\hbox{is}& \hbox{a commutative subfield of $D$},\\
S &=&D(x)\langle\p\rangle,\\
S^{(m)} &=& S\e_1+...+S\e_m \hbox{, a free $S$-module of rank $m$}, \\
\de_1,...,\de_m &\hbox{is}& \hbox{a finite set of $K$-linearly
     independent elements in $K\langle x,\p\rangle$}, \\
\sigma(\alpha,f)&=&\sum_{i=1}^m\alpha\de_if\e_i \in S^{(m)},
     \ (\alpha\in S,f\in K\langle x,\p \rangle), \\
P(\alpha,f)&=&S\sigma(\alpha,f), \hbox{ ideal of } S^{(m)}, \\
{\cal D}_r &\hbox{is}&
     \hbox{the quotient ring of $A_r$}, \\
{\cal R}_r &=& {\cal D}_r (x_{r+1},...,x_n) \langle\p_{r+1},...,\p_n\rangle, \\
{\cal S}_r &=& {\cal D}_r(x_{r+1},...,x_n)\langle\p_{r+1}\rangle.
\end{eqnarray*}
With exception of some minor changes we tried to stick to the
notation in \cite{bjork}.

\section{Preliminaries}%%%%%%%%%%%%%%%%%%%%%%%%%%%%%%%%%%%%%%%%%%%%%%%%%%%%%
Several useful properties of Weyl algebras are discussed in this
section. Also, we introduce a few rings that will come handy later
on.
\subsection{$A_n$ is simple}%%%%%%%%%%%%%%%%%%%%%%%%%%%%%%%%%%%%%%%%%%%%%%%%
To see that $A_n$ is simple, i.e. has no nontrivial two-sided
ideals, we notice that, for $f=\sum_i x^{\alpha_i}\p^{\beta_i}\in
A_n\setminus\{0\}$ in the standard form, ${df}/{dx_r} = \p_r f -
f\p_r$ for $r=1,...,n$, where ${\p f}/{\p x_r}$ is the formal
derivative of the above expression of $f$ with respect to $x_r$.
Similarly, ${df}/{d\p_r} = fx_r - x_rf$ for the formal derivative
with respect to $\p_r$. Note that these formal derivatives as well as all
the multiple derivatives of $f$ belong
to the two-sided ideal $A_nfA_n$.

Now assume $x^\alpha\p^\beta$ is the leading term of $f$ with
respect to some total degree monomial ordering. We are going to
perform $|\alpha|+|\beta|$ differentiations: for all $i=1,...,n$
differentiate $f$ $\alpha_i$ times with respect to $x_i$ and
$\beta_i$ times with respect to $\p_i$. Under such operation the
leading term becomes equal to $\prod_{i=1}^n {\alpha_i!\beta_i!}$
and all the other terms vanish. Since the derivatives of $f$ don't
leave $A_nfA_n$, we showed that there is a simple algorithm to
find such $s_i,r_i\in A_n$ that
$$
\sum_{i=1}^m s_{i}fr_{i}=1.
$$
Hence, $A_nfA_n=A_n$, so $A_n$ is simple.

\subsection{$A_n$ is an Ore domain}%%%%%%%%%%%%%%%%%%%%%%%%%%%%%%%%%%%%%%%%%
\begin{proposition}
$A_n$ is an Ore domain, i.e. $A_nf\cap A_ng \neq 0$ and $fA_n\cap
gA_n \neq 0$ for every $f,g\in A_n\setminus\{0\}$.
\end{proposition}
\begin{proof}
See the proof of Proposition 8.4 in Bj\"ork \cite{bjork}.
\end{proof}

Let us point out that using Gr\"obner bases methods (see next
subsection) we can find a left(right) common multiple of $f,g\in
A_n\setminus\{0\}$, in other words we can find a nontrivial
solution to the equations $af=bg$ and $fa=gb$ where $a$ and $b$
are unknowns.

\subsection{Gr\"obner bases in $A_n$}
As we mentioned before, the notion of Gr\"obner basis of a (left)
ideal can be defined for Weyl algebras in the same way as it is
defined in the case of polynomials. Moreover, Buchberger algorithm
for computing Gr\"obner bases works, leading to algorithms for
computing intersections of ideals, kernels of maps, syzygy
modules, etc. A good reference on Gr\"obner bases for algebras of
solvable type is \cite{nonCommGrobner}.

\subsection{More rings}
There is a quotient ring $D$ associated to every Ore domain $A$. Ring $D$
is a skew field that can be constructed both as the ring of left fractions
$a^{-1}b$ and as the ring of right fractions $cd^{-1}$,
where $a,b,c,d \in A$. There is a detailed treatment of this issue in
\cite{bjork}.

Let $D$ be a skew field, we will be interested in the ring
$S=D(x)\langle\p\rangle$, which is a ring of differential operators
with coefficients in $D(x)$. It is easy to see that $S$ is simple.

Since the Weyl algebra $A_r$ is an Ore domain, we can form its
quotient ring, which we denote by ${\cal D}_r$. The $S$ we are
going to play with is ${\cal S}_r = {\cal
D}_r(x_{r+2},...,x_n)(x_{r+1})\langle\p_{r+1}\rangle$. Let us
state without proof a proposition which shall help us to compute
Gr\"obner bases in ${\cal S}_r$.
\begin{proposition}
Let $F=\{f_1,...,f_k\}\subset A_n$ is a generating set of left ideal $I$ of
${\cal S}_r$. Compute a Gr\"obner basis $G=\{g_1,...,g_m\}$ of $A_n\cdot F$
with respect to any monomial ordering eliminating $\p_{r+1}$. Then $G$ is
contained in ${\cal S}_r \cap A_n$ and is a Gr\"obner basis of $I$.
\end{proposition}

\section{Holonomic modules are cyclic}%%%%%%%%%%%%%%%%%%%%%%%%%%%%%%%%%%%%%%
In this section we consider a simple ring $A$ such that $A$ has finite
length as a left module over itself.
Note that $A_n$ is such a ring.

\begin{theorem}\label{tOne}
Every left $A$-module $M$ of finite length is cyclic.
In particular every holonomic $A_n$ module is cyclic.
\end{theorem}

Suppose we know how to compute a cyclic generator
for every module $M'$ of length
less than $l$. For length $0$ such generator would be $0$.

Consider a module $M$ of length $l$. Take $0 \neq \alpha\in M$.
If $M=A\alpha $ then we are done. If not then since $l(M/A\alpha )<l$  by
induction we can find $\beta $ such that its image in $M/A\alpha $ is a
cyclic generator. Now $M = A\cdot \{\alpha,\beta\}$ and what we need to prove
is
\begin{lemma} \label{lCyclic}
Let $M$ be a left $A$-module of finite length
and $\alpha ,\beta \in M$. Then there
exists $\gamma \in M$ such that $A\gamma =A\alpha +A\beta $.
\end{lemma}
\begin{proof}
Define two functions $l_1$ and $l_2$ for pair $(\alpha,\beta)$.
\begin{eqnarray*}
l_{1}(\alpha,\beta) &=& \len(A\beta) \\
l_{2}(\alpha ,\beta ) &=&
\len\left(\left( A\alpha +A\beta \right) /A\alpha\right).
\end{eqnarray*}

Let also introduce an order $<$ on the set of pairs
$(\alpha,\beta)\in M\times M$:
\begin{eqnarray*}
(\alpha',\beta')<(\alpha,\beta) &\Leftrightarrow&
(l_1(\alpha',\beta'),l_2(\alpha',\beta'))
<_{lex}
(l_1(\alpha,\beta),l_2(\alpha,\beta))\\
&\Leftrightarrow&
l_1(\alpha',\beta')<l_1(\alpha,\beta) \\
&&\hbox{ OR } \left( l_1(\alpha',\beta')=l_1(\alpha,\beta) \hbox{ AND }
l_2(\alpha',\beta')<l_2(\alpha,\beta)
\right)
\end{eqnarray*}

Suppose for any pair $(\alpha',\beta') < (\alpha,\beta)$, we can find
$\gamma'\in M$ such that $A\gamma' = A\cdot \{\alpha',\beta'\}$.

Let the ideals $L(\alpha )$ and $L(\beta )$ in $A$ be the annihilators of $%
\alpha $ and $\beta $ respectively.
Since $\len(A) = \infty$, we know that $L(\alpha )\neq 0$;
pick any element $0\neq f\in L(\alpha )$. Since $A$ is simple we can find
$s_{i},r_{i}\in A$, $I=1,...,M$ such that
\begin{equation}
\sum_{i=1}^m s_{i}fr_{i}=1. \label{eSR}
\end{equation}

Consider two cases:

\begin{enumerate}
\item  There is some $r=r_{i}$ such that $L(\beta )+L(\alpha )r=A.$
\item  The opposite is true.
\end{enumerate}

Case 1.
We can write $1=E_{\alpha }r+E_{\beta }$ for some $E_{\alpha
},E_{\beta }\in A$ such that $E_{\alpha }\alpha =0$ and $E_{\beta }\beta =0$%
. Let $\gamma =\alpha +r\beta $.

Now we can get $\beta $ from $\gamma $:
\begin{equation*}
\beta =\left( E_{\alpha }r+E_{\beta }\right) \beta =E_{\alpha }r\beta
=E_{\alpha }\alpha +E_{\alpha }r\beta =E_{\alpha }\gamma \text{.}
\end{equation*}
Hence $\beta \in A\gamma $ and since $\alpha =\gamma -r\beta $ the module $%
M=A\alpha +A\beta $ is indeed generated by $\gamma $.

Case 2.
From (\ref{eSR}) it follows that $\sum L(\beta )+Afr_{i}=A$, hence,
$\sum A(fr_{i}\beta )=A\beta $, so there is $r=r_{i}$
such that
\begin{equation}
A(fr\beta)\nsubseteq A\alpha. \label{eNotInAlpha}
\end{equation}

Since we are not in case 1,
$
L(\beta )+Afr\subset L(\beta )+L(\alpha )r\neq A
$.
Take this modulo $L(\beta )$ to get
\begin{equation}
A(fr\beta ) \cong \left( L(\beta )+Afr\right) /L(\beta )
\subsetneq A/L(\beta )\cong A\beta ,
\label{eProperBeta}
\end{equation}
so $A(fr\beta )$ is proper in $A\beta $.

The last statement implies $l_{1}(\alpha,fr\beta )<l_{1}(\alpha,\beta) $,
hence,
$(\alpha,fr\beta)<(\alpha,\beta)$,
so by induction hypothesis we can find $\gamma'\in M$
such that $A\gamma'=A(fr\beta )+A\alpha $.

Now (\ref{eNotInAlpha}) guarantees that  $l_{2}(\gamma',\beta
)<l_{2}(\alpha,\beta)$, and by induction we can find $\gamma$ for which
\begin{equation*}
A\gamma =A\gamma'+A\beta =A(fr\beta )+A\alpha +A\beta =A\alpha
+A\beta.
\end{equation*}
\end{proof}

\begin{remark}
There is an algorithm that finds a cyclic generator for a
holonomic left module over a Weyl algebra, since every step in the
proof of the Lemma \ref{lCyclic} is computable. The most
non-trivial and time consuming operation is producing the
annihilators $L(\alpha+r\beta)$ and $L(fr\beta)$ in the proof of
Lemma \ref{lCyclic} provided $L(\alpha)$ and $L(\beta)$. This
is done using Gr\"obner bases technique.
\end{remark}

We have programmed the algorithm corresponding to the proof of Theorem
\ref{tOne} using {\em Macaulay 2}.

\noindent {\bf Example.} Let us view the ring of polynomials
$k[x]$ as an $A_1$-module under the natural action of differential
operators. It has an irreducible module, because starting with a
nonzero polynomial $f$ we can obtain a nonzero constant by
differentiating it $\deg(f)$ times. The module $M=k[x]^3$ is the
direct sum of 3 copies of $k[x]$, is holonomic ($\len(M)=3$) and
is generated by vectors $(1,0,0),(0,1,0),(0,0,1)$. Our algorithm
produces a cyclic generator $\gamma=(x^2,x,1)$ and its
$A_1$-annihilator $L(\gamma)=A_1\partial^3$.

\section{Ideals are 2-generated}%%%%%%%%%%%%%%%%%%%%%%%%%%%%%%%%%%%%%%%%%%%%
In this section we give an effective proof of
\begin{theorem}\label{tTwo}
Every left ideal of the Weyl algebra $A_n$
can be generated by two elements.
\end{theorem}
\begin{proof}[Proof for $A_1$]
In this case the theorem follows from the fact that module $A_1/J$
is holonomic for any nonzero ideal $J$ of $A_1$.

Indeed, let $I$ be a left ideal of $A_1$. Pick $f\in I$ and set
$J=A_1f$. Then $I/J$ is a submodule of the holonomic module
$A_1/J$, hence, is holonomic. By Theorem \ref{tOne} there is
${\bar g} \in I/J$ such that $A_1{\bar g}=I/J$. Find a lifting
$g\in A_1$ such that $\bar g = g \mod J$. Elements $f$ and $g$
generate $I$.
\end{proof}

However, the theorem for $n>1$ makes a much tougher challenge.

\subsection{Lemmas for $S$} %%%%%%%%%%%%%%%%%%%%%%%%%%%
Let us explore some properties of $S=D(x)\langle\p\rangle$, the ring of linear
differential operators with coefficients in rational expressions in $x$ over
a skew field $D$.

Let $K$ be a commutative subfield of $D$, let $\de_1,...,\de_m$
be a finite set of
$K$-linearly independent elements in $K\langle x,\p\rangle \subset S$,
and let $S^{(m)} = S\e_1+...+S\e_m$ be a free $S$-module of rank $m$.

Also define $\sigma(\alpha,f)\in S^{(m)}$ to be the following sum
$\sigma(\alpha,f)=\sum_{i=1}^m\alpha\de_if\e_i,$ and
$P(\alpha,f)=S\sigma(\alpha,f)$ the submodule of $S^{(m)}$ generated
by $\sigma(\alpha,f)$. Note that $\sigma(\alpha,f)$ is $S$-linear
in $\alpha$ and respects addition in $f$.

\begin{lemma} \label{lSM}
Let $0\neq \alpha\in S$ and let $M$ be an $S$-submodule of
$S^{(m)}$ generated by $\{ \sigma(\alpha,f)| f\in K\langle
x,\p\rangle\}$. Then $M=S^{(m)}$.
\end{lemma}
\begin{proof}
%\iffalse
Without loss of generality let us assume that $\alpha\in D\langle
x,\p\rangle$: if not we can always find such $p\in D[x]$ that
$p\alpha\in D\langle x,\p\rangle$.

Fix a monomial ordering that respects the total degree in $x$ and
$\p$. For vector $v=\sum v_i\e_i \in (D\langle x,\p\rangle)^{(m)}$
denote by $\hbox{lm}(v)$ the largest of the the leading monomials
of the components $v_i$ of $v$ in this ordering.

Now start with vector $v=v^{(0)}=\sigma(\alpha,1)$; its components
$v_i=\alpha\de_i$ are $D$-linearly independent. Note that
computing expressions $\pi(v)=\p v-v\p$ and $\chi(v)=vx - xv$ has
an effect of differentiating each component of $v$ formally with
respect to $x$ and $\p$ respectively. These operations lower the
total degree of $v$ by $1$ if the differentiation is done with
respect to a variable that is present in $\hbox{lm}(v)$. Also, it
is not hard to see that they keep us in module $M$; for example,
for $v_0$ we have $\pi(v_0)=\p
v_0-v_0\p=\p\sigma(\alpha,1)-\sigma(\alpha,\p)$.

Run the following algorithm: initialize $v:=v_0$, while
$\hbox{lm}(v)$ contains an $x$ set $v:=\pi(v)$, then while
$\hbox{lm}(v)$ contains a $\p$ we set $v:=\chi(v)$. Since each
step lowers the total degree of $v$ by $1$, this procedure
terminates producing vector $w\in M$ of total degree $0$.

Hence, $w=w_{i_1}\e_{i_1}+...+w_{i_t}\e_{i_t}$ where $0\neq
w_{i_j}\in D$ for $j=1,...,t$. Via multiplying on the left by the
inverse of $w_{i_1}$ we can get the relation
\begin{equation}
 \e_{i_1}=a_2\e_{i_2}+...+a_t\e_{i_t} \label{lE}
\end{equation}
with $a_j\in D$ for $j=2,...,t$.

Now take $v^{(0)}$ and reduce it using (\ref{lE}). We get vector
$v^{(1)}$ whose $i_1$-th component is $0$ and the remaining
components are $D$-linearly independent, since the components of
$v^{(0)}$ are.

Repeat the above algorithm for $v=v^{(1)}$ and so on.
At the end we get a vector which is a scalar multiple of $\e_i$ for some $i$,
hence $e_i\in M$.
Using relations (\ref{lE}) we see that all basis vectors $\e_j$,
for $j=1,...,m$, are in $M$.
%\fi
\end{proof}

\begin{remark}\label{rFindF}
From the proof it follows that given a submodule $M$ of $S^{(m)}$ and
$\alpha\in S$ one can find $f\in K\langle x,\p\rangle$ such that
$\sigma(\alpha,f)\notin M$ algorithmically.
\end{remark}

The next lemma is central in the proof of the result. Note that
every step of the proof of the lemma can be carried out
algorithmically.
\begin{lemma}\label{lCentral}
Let $M$ be an $S$-submodule of $S^{(m)} = S\e_1+...+S\e_m$ such that
$\len(S^{(m)}/M)<\infty$. We can find $f\in K\langle x,\p\rangle$
such that $S^{(m)} = M + P(\alpha,f)$.
\end{lemma}
\begin{proof}
Let $l=\len(S^{(m)}/M)$. Assume the assertion is proved for all
$M'$ such that $\len(S^{(m)}/M')<l$. Remark \ref{rFindF} says that
we can find an $f\in K\langle x,\p\rangle$ such that
$\sigma(\alpha,f)$ doesn't belong to $M$.

For $t\in S$, $g\in K\langle x,\p\rangle$ let us define two $S$-modules
\begin{eqnarray*}
N_1 &=& M + P_1 \mbox{, where }P_1=P(\alpha,g), \\
N_2 &=& M + P_2 \mbox{, where }P_2=P(t\alpha,g).
\end{eqnarray*}

{\bf Claim.} {\em There is a module $M'$ such that
$M\subset M' \subset M+P(\alpha,f)$, $t\in S$,
and $g\in K\langle x,\p\rangle$ for which
\begin{eqnarray*}
&&t\sigma(\alpha,f) \in M, \\
&&M'+P(t\alpha,g)=S^{(m)}, \\
&&N_1=N_2.
\end{eqnarray*}
}

To prove this we employ (second) induction on $\len(M'/M)$. We
start with $M' = M+P(\alpha,f)$. We can find $0\neq t\in S$ such
that $t\alpha\sum{\de_if\e_i} \in M$; it follows from $S$ being
Ore. By the first induction hypothesis, for $M'$ and $t\alpha$
there exists $g\in K\langle x,\p\rangle$ such that $M' +
P(t\alpha,g) = S^{(m)}$. Notice that $N_1\supset N_2$ and $M'+P_i=
S^{(m)}$ for $i=1,2$. Also for $i=1,2$ we have
    $$S^{(m)}/N_i = (M'+P_i)/(M+P_i) = M'/(M+M'\cap P_i).$$ \label{eN1N2}
If $\len(S^{(m)}/N_1)=\len(S^{(m)}/N_2)$ then $N_1=N_2$ and we are
done. We are done as well if $N_1=S^{(m)}$. If both conditions
above fail, by looking at the right hand side of \ref{eN1N2} we
determine that $M''=M+M'\cap P_1$ both contains $M$ and is
contained in $M'$ properly, plus $\len(M''/M)<\len(M'/M)$. Set
$M':=M''$ and repeat the above procedure.

\medskip
To finish the proof of the lemma we take $M',t,g$ as in the claim and assert
that $N' = M + P(\alpha,f+g)$ equals $S^{(m)}$.
Indeed, $\sigma(t\alpha,f+g)=t\sigma(\alpha,f)+\sigma(t\alpha,g)=
\sigma(t\alpha,g)$ modulo $M$, so $N_2\subset N'$. But $N_1=N_2$,
thus $\sigma(\alpha,g)\in N'$, hence,
$\sigma(\alpha,f) = \sigma(\alpha,f+g) - \sigma(\alpha,g) \in N'$.
Now we see that $M'\subset N'$ and $P_2\subset N'$. Since $M'+P_2=S^{(m)}$,
we proved $N'=S^{(m)}$.
\end{proof}

\subsection{Lemmas for ${\cal R}_r$} %%%%%%%%%%%%%%%%%%%%%%%%%%%%%%%%%%%%
At this stage we shall specify the components in the definition of
$S=D(x)\langle\p\rangle$. We set $D={\cal D}_r(x_{r+2},...,x_n)$,
$x=x_{r+1}$ and $\p=\p_{r+1}$, so that new $S$ is equal to ${\cal
S}_r={\cal D}_r(x_{r+1},x_{r+2},...,x_n)\langle \p_{r+1} \rangle $
which is a subring of ${\cal R}_r$. Also the commutative subfield
$K$ of $D$ that showed up before is replaced by the $k$, the
coefficient field from the definition of $A_n=A_n(k)$.

\begin{proposition} \label{pMain}
Let $\de_1,...,\de_m$ be a finite set of $K$-linearly independent
elements in $K\langle x_{r+1},\p_{r+1}\rangle$
and let $0 \neq \rho \in A_{r+1}[x_{r+2},...,x_n]$.
Let $S^{(m+1)} = S\e_0+S\e_1+...+S\e_m$ be a free $S$-module of rank $m+1$
And let $S^{(m+1)}\rho \subset S^{(m+1)}$ be its $S$-submodule generated by
$\{\rho \e_0, \rho \e_1,...,\rho \e_2\}$. Then there exists some $f\in K $ such that
    $$ S^{(m+1)} = S^{(m+1)}\rho + S(\e_0+\de_1f\e_1+...+\de_mf\e_m). $$
\end{proposition}
\begin{proof}
Follows from Lemma \ref{lCentral}
\end{proof}

\begin{lemma} \label{l85}
Let $q\in A_r[x_{r+1},...,x_n]$
and let $a_1,...,a_t$ be a finite set in $A_n$.

Then there exists some $0\neq\rho\in A_r[x_{r+1},...,x_n]$
such that $\rho a_j \in A_n q$ for all $j$.
\end{lemma}
\begin{proof}
See the proof of Lemma 8.5 in Bj\"ork \cite{bjork}.
\end{proof}

Let us point out that once we know the statement of the lemma is true,
we can compute the required $\rho$ by finding a Gr\"obner basis of the
module of syzygies of the columns of the matrix
$$ \left( \begin{array}{ccccc}
a_1 & q & 0 &...& 0 \\
a_2 & 0 & q &...& 0 \\
... &...&...&...&...\\
a_t & 0 & 0 &...& q
\end{array} \right) $$
with respect to a monomial order that eliminates
$\p_{r+1},...,\p_n$ and such that $\e_1>\e_2>...>\e_{t+1}$ where
$\e_1,\e_2,...,\e_{t+1}$ is the basis ($\e_i$ corresponds to the
$i$-th column) of the free module $A_n^{t+1}$ containing our
submodule of syzygies. Such Gr\"obner basis is guaranteed (by
Lemma \ref{l85}) to contain some syzygy producing the relation
$\rho\e_1+b_2\e_2+...+b_{t+1}\e_{t+1}=0$ where $\rho \in
A_r[x_{r+1},...,x_n]$, $b_i\in A_n$ for $i=2,...,n$. It is not
hard to see that this is the $\rho$ we need.

\begin{lemma}
Let $0 \neq q \in A_{r+1}[x_{r+2},...,x_n]$ and
let $u,v\in A_n$ with $v \neq 0$.
Then there is some $f\in A_n$ such that
${\cal R}_r = {\cal R}_r q + {\cal R}_r(u+vf)$.
\end{lemma}
\begin{proof}
Consider  the following subring of $A_n$ obtained by "removing" $x_{r+1}$
and $\p_{r+1}$:
$$ A_{\widehat{r+1}} =
k\langle x_1,...x_r,x_{r+2},...,x_n,\p_1,...,\p_r,\p_r+2,...,\p_n\rangle.$$
Now $A_n = A_{\widehat{r+1}} \otimes_k k\langle x_{r+1},\p_{r+1}\rangle$, so we can write $v = \de_1g_1+...+\de_mg_m$
where $\de_1,...,\de_m$ are elements of $k\langle x_{r+1},\p_{r+1}\rangle$
linearly independent over $k$ and $g_1,...,g_m\in A_{\widehat{r+1}}$.
The ring $A_{\widehat{r+1}}$ is simple, since it is a Weyl algebra,
thus we can find such $h_1,...,h_l\in A_{\widehat{r+1}}$ that
$$A_{\widehat{r+1}} = \sum_{i=1}^m\sum_{j=0}^l A_{\widehat{r+1}}g_ih_j.$$
Since $A_{\widehat{r+1}}$ is a subring of ${\cal R}_r$ it means that
${\cal R}_r = \sum\sum {\cal R}_r g_ih_j$.

{\bf Sublemma.} {
\em For any $b_1,...,b_m \in A_{\widehat{r+1}}$ there exists some
$f\in k\langle x_{r+1},\p_{r+1}\rangle$ such that
$$ {\cal R}_r q + {\cal R}_r u + {\cal R}_r b_1 + ... + {\cal R}_r b_m
= {\cal R}_r q + {\cal R}_r (u+\de_1 f b_1 + ... + \de_m f b_m).$$
}
\begin{proof}
It follows from Lemma \ref{l85} that there is
$0\neq\rho\in A_r[x_{r+1},...,x_n]$
such that $\rho b_1,...,\rho b_m \in A_n q$ as well as $\rho u \in A_n q$.
With the help from Proposition \ref{pMain} we get
$f\in k\langle x_{r+1},\p_{r+1}\rangle$ such that
$ S^{(m+1)} = S^{(m+1)}\rho + S(\e_0+\de_1f\e_1+...+\de_mf\e_m) $
and since $S$ is a subring of ${\cal R}_r$ we have
\begin{equation}
{\cal R}_r^{(m+1)} = {\cal R}_r^{(m+1)}\rho +
    {\cal R}_r(\e_0+\de_1f\e_1+...+\de_mf\e_m). \label{eUMM}
\end{equation}
Now map $\e_0 \mapsto u$ and $\e_i \mapsto b_i$ for all $i$; this
map from ${\cal R}_r^m$ to ${\cal R}_r$ has its image equal to $
{\cal R}_r q + {\cal R}_r u + {\cal R}_r b_1 + ... + {\cal R}_r
b_m$ and maps the right hand side of (\ref{eUMM}) to a subset of
${\cal R}_r q + {\cal R}_r (u+\de_1 f b_1 + ... + \de_m f b_m)$,
because $\rho u, \rho b_1, ..., \rho b_m \in A_n q$. Moreover
these two expressions are equal, since it is easy to see that the
latter is contained in the former as well.
\end{proof}

{\em Proof of lemma continued.} We apply our Sublemma to
$b_i=g_ih_1$ ($i=1,...,m$)
to get $f_1\in k\langle x_{r+1},\p_{r+1}\rangle$ such that
$$
{\cal R}_r q + {\cal R}_r u + \sum_{j=1}^m {\cal R}_r g_ih_1
= {\cal R}_r q + {\cal R}_r (u+\sum_{j=1}^m\de_i f_1 g_ih_1).
$$
Since $v=\de_1g_1+...+\de_mg_m$ and since $f_1$ commutes with all $g_i$,
the last equation transforms into
$$
{\cal R}_r q + {\cal R}_r u + \sum {\cal R}_r g_ih_1
= {\cal R}_r q + {\cal R}_r (u+v f_1 h_1).
$$

Now reapply the Sublemma with $u$ replaced by $u+v f_1 h_1$ and
$b_i=g_ih_2$ ($i=1,...,m$). As in the first step we get
\begin{eqnarray*}
&&{\cal R}_r q + {\cal R}_r u + \sum {\cal R}_r g_ih_1 + \sum {\cal R}_r g_ih_2 \\
&&= {\cal R}_r q + {\cal R}_r (u+v f_1 h_1) + \sum {\cal R}_r g_ih_2 \\
&&= {\cal R}_r q + {\cal R}_r (u+v f_1 h_1+f_2 h_2)
\end{eqnarray*}
for some $f_2\in k\langle x_{r+1},\p_{r+1}\rangle$. After $l$ many steps we
arrive at
$$
{\cal R}_r
= {\cal R}_r q + {\cal R}_r u + \sum_{i=1}^m\sum_{j=1}^l {\cal R}_r g_ih_j
= {\cal R}_r q + {\cal R}_r (u+v \sum_{j=1}^l  f_i h_i),
$$
which proves the lemma with $f=f_1h_1+...+f_lh_l$.
\end{proof}

The following lemma follows from the previous one.
\begin{lemma} \label{lMain}
Let $0\leq r \leq n-1$ and let $0 \neq q \in A_{r+1}[x_{r+2},...,x_n]$
and let $u,v\in A_n$ with $v \neq 0$.
Then there is some $f\in A_n$,$q'\in A_r[x_{r+1},...,x_n]$ such that
$q' \in A_nq+A_n(u+vf)$.
\end{lemma}
\begin{proof}
It is easy to see that this lemma is equivalent to the previous one.
\end{proof}

\subsection{Final chords}%%%%%%%%%%%%%%%%%%%%%%%%%%%%%%%%%%%%%%%%%%%%%%%%%%%%

\begin{proposition}[r]
Let $0\leq r \leq n$, there is some $q_r\in A_r[x_{r+1},...,x_n]$
and $d_r,e_r\in A_n$ such that $q_r c \in A_n(a+d_r c)+A_n(b+e_r
c)$.
\end{proposition}
\begin{proof}
The statement is true for $r=n$, since $A_n$ is Ore and
$A_nc\cap(A_na+A_nb)$.

Fix $r$. Assume that the statement is true for $r+1,...,n$, then
there exist $q_{r+1},d_{r+1},e_{r+1}$ such that $q_{r+1} c \in A_n
a'+A_n b'$, where $a' = a+d_{r+1} c$ and $b' = b+e_{r+1} c$. Hence
we can write $q_{r+1}c=h_1 a' + h_2 b'$, where we can take $h_1
h_2\neq 0$ since $A_n a' \cap A_n b' \neq 0$. Also since $h_1 A_n
\cap h_2 A_n \neq 0$ we can also find $g_1,g_2$ satisfying
$h_1g_1+h_2g_2=0$, and since $A_n q_{r+1}c \cap A_n b'\neq 0$
there are $s,t$ such that $sq_{r+1}c = tb'$ . Using Lemma
\ref{lMain} to $q=q_{r+1}$ with $u=0$ and $v=tg_2$, we get
$q_r=q'$ and $f$ such that $q_r = p_1q_{r+1} + p_2tg_2f$ for some
$p_1,p_2$. Summarizing, there exist such
$h_1,h_2,g_1,g_2,s,t,p_1,p_2 \in A_n\setminus\{0\}$ that
\begin{eqnarray*}
&&q_r = p_1q_{r+1} + p_2tg_2f \\
&&q_{r+1}c=h_1 a' + h_2 b' \\
&&h_1g_1+h_2g_2=0 \\
&&sq_{r+1}c = tb'
\end{eqnarray*}

Using these 4 equations, make the following calculation: (In each
section the underlined terms sum up to $0$.)
\begin{eqnarray*}
q_r c & = &p_1q_{r+1}c + p_2tg_2fc \\
      & = &p_1q_{r+1}c - \underline{p_2sq_{r+1}c} \\
      & + &p_2tg_2fc + \underline{p_2tb'} \\
      & = &(p_1-p_2s)q_{r+1}c + p_2t(b'+g_2fc)\\
      & = &(p_1-p_2s)(h_1 a' + h_2 b') + p_2t(b'+g_2fc)\\
      & = &(p_1-p_2s)h_1 a' + \underline{(p_1-p_2s)h_1g_1fc}\\
      & + &(p_1-p_2s)h_2b' + \underline{(p_1-p_2s)h_2g_2fc} + p_2t(b'+g_2fc)\\
      & = &(p_1-p_2 s)h_1(a' + g_1f c) + ((p_1-p_2 s)h_2+p_2t)(b' + g_2f c).
\end{eqnarray*}
Thus, with $d_r=d_{r+1} + g_1f c$ and $e_r=e_{r+1} + g_2f c$ the
conclusion of the proposition holds.
\end{proof}

The proposition above (for $r=0$) shows that by ``elimination'' of
variables $\p_i$ one at a time we can get such $d,e\in A_n$ that
$q_0c \in A_n(a+dc)+A_n(b+ec)$ where $q_0\in k[x_1,...,x_n]$. This
proves a 50\% version of Theorem \ref{tTwo}:
\begin{theorem}\label{t50}
Every ideal of $k(x_1,...,x_n)\langle \p_1,...,\p_n \rangle$
can be generated by two elements.
\end{theorem}

To go other 50\% of the way one has to do a similar kind of
``elimination'' of $x_i$-s. This amounts to making copies of all
lemmas that we stated for a slightly different set of rings. The
trickiest part is considering ring ${{\cal
S}_r}'=k(x_1,...,x_r)\langle x_{r+1},\p_{r+1} \rangle$ instead of
${\cal S}_r$. In other words instead of a ring of type
$D(x)\langle\p\rangle$ where $D$ is a skew field, we have to
consider the first Weyl algebra $A_1({\cal K})$ where ${\cal K}$
is a (commutative) field. Fortunately, analogues of Lemmas
\ref{lSM} and \ref{lCentral} for the latter ring can be
effectively proved along the same lines.

\noindent {\bf Examples.} (1) Consider $A_3$. For
$a=\p_1$,$b=\p_2$,$c=\p_3$ one can show that $A_3\cdot \{a,b,c\} =
A_3\cdot \{a,b+x_1c\}$. Indeed, the following calculation displays
it:
$$ c = (-x_1\p_3-\p_2)a + \p_1(b+x_1c).$$

(2) Another example is produced by our algorithm implemented in
{\em Macaulay 2}. Let $a=\p_1+x_3$, $b=\p_2^2+x_2+x_3^2$,
$c=\p_3+x_1$. Then the ideal $A_3\cdot\{a,b,c\}$ is generated by
$\p_1+x_3$ and
$\p_2^2+(x_1^2x_3+x_1)\p_3+x_1^3x_3+x_1^2+x_3^2+x_2$.

(3) In case of $A_1$ we can construct a more efficient algorithm
based on the proof of Theorem \ref{tTwo} given for this special
case. Here is a {\em Macaulay 2} script computing 2 generators for
the annihilating ideal $I \subset A_1(\mathbb{Q})$ of the set of
polynomials $\{ax^4+bx^6+cx^8+dx^{10}\ |\ a,b,c,d\in
\mathbb{Q}\}\subset \mathbb{Q}[x]$.
\begin{Macaulay2}
\begin{verbatim}
i1 : load "D-modules.m2"; load "stafford.m2";

i3 : R = QQ[x,D, WeylAlgebra=>{x=>D}];

i4 : L = {4,6,8,10};

i5 : I = ideal gens gb intersect apply(L, i->PolyAnn x^i);

             4 4      3 3       2 2                   11
o5 = ideal (x D  - 22x D  + 207x D  - 975x*D + 1920, D  , ...

o5 : Ideal of R

i6 : time J = ideal stafford I
     -- used 73.08 seconds

             4 4      3 3       2 2
o6 = ideal (x D  - 22x D  + 207x D  - 975x*D + 1920,

             3 15    2 15      2 14       14        13
            x D   + x D   + 15x D   + 9x*D   + 58x*D

                       13      12    11
                  + 15D   + 50D   + D  )

o6 : Ideal of R

i7 : I == J

o7 = true
\end{verbatim}
\end{Macaulay2}

\section{Conclusion}
The implementations of
the algorithms constructed along the lines of the proofs of
Theorems \ref{tOne} and \ref{tTwo} in {\em Macaulay 2} work
only on rather small examples
for quite obvious reason: the expression swell in Gr\"obner
bases computations.

Let us comment on the differences of algorithm of Hillebrand and
Schmale \cite{schmale} and ours. Their algorithm takes care of
(weaker) Theorem \ref{t50}. As a step it includes enumerating a
certain infinite subset of polynomials in one variable and testing
them to satisfy a certain property, where the testing procedure
involves Gr\"obner bases computations. Although we
believe that their argument could be extended to build an
algorithm for 100\% of Stafford's theorem, it looks as the ``test
set'' for the remaining 50\% will be significantly more
complicated. Hence, our constructive approach at every step of the
algorithm seems to be more practical. Having programmed Hillebrand
and Schmale's algorithm as well, we have to point out, that it
faces the same type of expression swell as our program, hence
the comparison of performance is just a theoretical question at this
point.

\smallskip
Finally, let us mention that the algorithm for finding a cyclic
generator of a holonomic module is already included in the
\emph{D-modules package for Macaulay 2} \cite{Dmodules};
eventually, the algorithms for finding two generators of a
$A_n$-ideal will be added to the package as well.

\bibliographystyle{plain}

\begin{thebibliography}{1}

\bibitem{bjork}
Bj\"ork, J.-E.
{\em Rings of differential operators.}
 North-Holland Mathematical Library, 21. North-Holland Publishing Co.,
Amsterdam-New York, 1979.

\bibitem{MTwo}
Grayson, Daniel; Stillman, Michael.
{\em Computer algebra system Macaulay 2. }
{\tt http://www.math.uiuc.edu/Macaulay2}

\bibitem{schmale}
Hillebrand, Andre; Schmale, Wiland.
{\em Towards an effective version of a theorem of Stafford. }
Effective methods in rings of differential operators.
J. Symbolic Comput. 32 (2001), no. 6, 699--716.

\bibitem{nonCommGrobner}
Kandri-Rody, A.; Weispfenning, V.
{\em Noncommutative Gr\"obner bases in algebras of solvable type.}
J. Symbolic Comput. 9 (1990), no. 1, 1--26.

\bibitem{Dmodules}
Leykin, Anton; Tsai, Harrison.
{\em $D$-modules for Macaulay 2.}
\newline {\tt http://www.math.umn.edu/\~{}leykin/Dmodules}

\bibitem{stafford}
Stafford, J. T.
{\em Module structure of Weyl algebras.}
J. London Math. Soc. (2) 18 (1978), no. 3, 429--442.

\end{thebibliography}

\end{document}